\documentstyle{article}
\def\con(#1){\special{em:line #1}\kern0pt}
\def\po(#1,#2,#3){\put(#1,#2){\special{em:point #3}}}
\begin{document}
\hfill q-alg/9711024
\vspace{5mm}
\begin{center}
{ \large \bf FRT  quantization  theory  for the nonsemisimple
 Cayley-Klein groups}
\end{center}
\begin{center}
          N.A.Gromov, I.V.Kostyakov, V.V.Kuratov \\
    {\it  Department of Mathematics, Komi Science Center,}\\
   {\it Ural Division, Russian Academy of Science,}\\
  {\it 167610 Syktyvkar, RUSSIA}
\end{center}

{\bf Abstract.}
The quantization theory of the simple Lie groups and algebras was
developed by Faddeev-Reshetikhin-Takhtadjan (FRT).
In group theory there is a remarkable set of groups, namely
the motion groups of n-dimensional spaces of constant curvature
or the orthogonal Cayley-Klein (CK) groups.
 In some sense the CK
groups are in the nearest neighborhood with the simple ones.
The well known groups of physical
interest such as
Euclidean $ E(n), $ Poincare $ P(n), $ Galileian $ G(n) $ and
other nonsemisimple groups are in the set of CK groups.
   But many standart algebraical costructions are not suitable
for the nonsemisimple groups and algebras, in particular
Killing form is degenerate, Cartan matrix do not exist.
Nevertheless it is possible to describe and to quantize all
CK groups and algebras, as it was made for the simple ones.
   The principal proposal is to consider CK groups as the groups
over an associative algebra $ {\bf D} $ with nilpotent
commutative generators
and the corresponding quantum CK
groups as the algebra of noncommutative functions over $ {\bf D}. $

\section{Introduction}
It is well known [6], [8] that there are $ 3^n $
n--dimensional real spaces of constant curvature or
Cayley--Klein spaces. These spaces are the most symmetric ones,
i.e. their motion groups have the maximal dimension
and for this reason are often used in physics.
The simple group $ SO(n+1) $
is the motion group of the n--dimensional spherical space.
All other CK groups have the same dimension
$ n(n-1)/2 $ and may be obtained from $ SO(n+1) $
by the contractions and analytical continuations [2].

The notion of Lie group contraction was first introduced by
E.In\"on\"u and E.P.Wigner [9] as some limiting procedure
and was later extend on new algebraical structures
such as Lie bialgebra [18], Hopf algebra [11], graded contractions
[16], [17], but the fundamental idea of degenerate transformations is
presented in all cases (see [15] for detailes). On the
other hand the degenerate transformation is something
incorrect  from mathematical point of view. So
it is necessary to find instead of it an relevant
mathematical construction. It seems that the
consideration of the Lie groups as the groups over an
associative algebra {\bf D} with nilpotent
commutative generators is the appropriate tool  at least in CK
scheme. The validity of such approach is demonstrated
for  the FRT quantization theory [1]. It is possible  to
reformulate the
quantum deformations of the simple groups in such a way
to obtain the quantum deformations of all contracted
CK groups.

This paper is orginezed as follows. In Section 2 the
standart In\"on\"u--Wigner contraction and our approach are compared
for the simplest case of one dimensional CK spaces and their
motion group. The algebra  $ {\bf D}_n(\iota; {\bf C}) $
is introduced in Section 3. The orthogonal CK groups
$ SO(N;j;{\bf R}) $ are described in Section 4 as the matrix
groups over $ {\bf D}(\iota;{\bf R}). $
Section 5 is devoted to the quantum  orthogonal
CK groups $ SO_v(N; j; {\bf C}) $ and
the quantum algebras
$ so_v(N;j;{\bf C}) $ are obtained as the
dual object to the corresponding
quantum groups in Section 6. The developed
approach is illustrated in Sections 7 on the
example of $ N=3 $ quantum groups and algebras.
The final remarks are given in Conclusion.

\section{In\"on\"u-Wigner contractions: traditional approach
and suitable mathematical structure}

Let us regard finite rotations on the angle $ \varphi $
for three planes: euclidean (Fig.1), galileian (Fig.2) and
minkowskian (Fig.3).

\unitlength=0.50mm
\begin{figure}[hb]
\begin{picture}(210,90)(-80,0)
\put(83.00,54){\makebox(0,0)[cd]{$ x_1 $ }}
\put(57,84){\makebox(0,0)[cd]{$ x_0 $ }}
\put(50,50){\vector(1,0){35}}\put(50,50){\vector(0,1){35}}
\put(50,50){\line(-1,0){35}}\put(50,50){\line(0,-1){35}}
\put(50,50){\line(2,3){14}}
\put(56,64){\makebox(0,0)[cc]{$\phi$ }}
\po(50,75,1)\po(51,74.98,2)\po(52,74.92,3)\po(53,74.82,4)
\po(54,74.68,5)\po(55,74.5,6)\po(56,74.27,7)\po(57,74,8)
\po(58,73.69,9)\po(59,73.32,10)\po(60,72.91,11)\po(61,72.45,12)
\po(62,71.93,13)\po(63,71.35,14)\po(64,70.71,15)\po(65,70,16)
\po(66,69.21,17)\po(67,68.33,18)\po(68,67.35,19)\po(69,66.25,20)
\po(70,65,21)\po(71,63.56,22)\po(72,61.87,23)\po(73,59.8,24)
\po(74,57,25)\po(74.3,55.87,26)\po(74.6,54.45,27)\po(74.8,53.16,28)
\po(74.9,52.23,29)\po(75,50,30)
\po(50,25,59)\po(51,25.02,58)
\po(52,25.08,57)\po(53,25.18,56)\po(54,25.32,55)\po(55,25.5,54)
\po(56,25.73,53)\po(57,26,52)\po(58,26.31,51)\po(59,26.68,50)
\po(60,27.09,49)\po(61,27.55,48)\po(62,28.07,47)
\po(63,28.65,46)\po(64,29.29,45)\po(65,30,44)\po(66,30.79,43)
\po(67,31.67,42)\po(68,32.65,41)\po(69,33.75,40)\po(70,35,39)
\po(71,36.44,38)\po(72,38.13,37)\po(73,40.2,36)\po(74,43,35)
\po(74.3,44.13,34)\po(74.6,45.55,33)\po(74.8,46.84,32)
\po(74.9,47.77,31)
\po(49,25.02,60)\po(48,25.08,61)\po(47,25.18,62)
\po(46,25.32,63)\po(45,25.5,64)\po(44,25.73,65)
\po(43,26,66)\po(42,26.31,67)\po(41,26.68,68)
\po(40,27.09,69)\po(39,27.55,70)\po(38,28.07,71)
\po(37,28.65,72)\po(36,29.29,73)\po(35,30,74)
\po(34,30.79,75)\po(33,31.67,76)\po(32,32.65,77)
\po(31,33.75,78)\po(30,35,79)\po(29,36.44,80)
\po(28,38.13,81)\po(27,40.2,82)\po(26,43,83)
\po(25.7,44.13,84)\po(25.4,45.55,85)\po(25.2,46.84,86)
\po(25.1,47.77,87)
\po(25.7,55.87,92)\po(25.4,54.45,91)\po(25.2,53.16,90)
\po(25.1,52.23,89)
\po(49,74.98,116)\po(48,74.92,115)\po(47,74.82,114)
\po(46,74.68,113)\po(45,74.5,112)\po(44,74.27,111)
\po(43,74,110)\po(42,73.69,109)\po(41,73.32,108)
\po(40,72.91,107)\po(39,72.45,106)\po(38,71.93,105)
\po(37,71.35,104)\po(36,70.71,103)\po(35,70,102)
\po(34,69.21,101)\po(33,68.33,100)\po(32,67.35,99)
\po(31,66.25,98)\po(30,65,97)\po(29,63.56,96)
\po(28,61.87,95)\po(27,59.8,94)\po(26,57,93)\po(25,50,88)
\con(1,2)\con(2,3)\con(3,4)\con(4,5)
\con(5,6)\con(6,7)\con(7,8)\con(8,9)
\con(9,10)\con(10,11)\con(11,12)\con(12,13)
\con(13,14)\con(14,15)\con(15,16)\con(16,17)
\con(17,18)\con(18,19)\con(19,20)\con(20,21)
\con(21,22)\con(22,23)\con(23,24)\con(24,25)
\con(25,26)\con(26,27)\con(27,28)\con(28,29)
\con(29,30)\con(30,31)\con(31,32)\con(32,33)
\con(33,34)\con(34,35)\con(35,36)
\con(36,37)\con(37,38)\con(38,39)\con(39,40)
\con(40,41)\con(41,42)
\con(42,43)\con(43,44)\con(44,45)\con(45,46)
\con(46,47)\con(47,48)\con(48,49)\con(49,50)
\con(50,51)\con(51,52)\con(52,53)
\con(53,54)\con(54,55)\con(55,56)\con(56,57)\con(57,58)\con(58,59)
\con(59,60)\con(60,61)\con(61,62)\con(62,63)\con(63,64)\con(64,65)
\con(65,66)\con(66,67)\con(67,68)\con(68,69)\con(69,70)\con(70,71)
\con(71,72)\con(72,73)\con(73,74)\con(74,75)\con(75,76)\con(76,77)
\con(77,78)\con(78,79)\con(79,80)\con(80,81)
\con(81,82)\con(82,83)\con(83,84)\con(84,85)\con(85,86)\con(86,87)
\con(87,88)\con(88,89)\con(89,90)\con(90,91)
\con(91,92)\con(92,93)\con(93,94)\con(94,95)\con(95,96)\con(96,97)
\con(97,98)\con(98,99)\con(99,100)\con(100,101)
\con(101,102)\con(102,103)\con(103,104)
\con(104,105)\con(105,106)\con(106,107)
\con(107,108)\con(108,109)\con(109,110)\con(110,111)
\con(111,112)\con(112,113)\con(113,114)\con(114,115)\con(115,116)
\con(116,1)
\end{picture}
\caption{Rotations on the euclidean plane}
\end{figure}
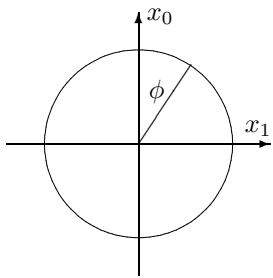

\begin{figure}[htp]
\begin{picture}(200,90)(-80,0)
\put(50,50){\line(1,2){12.5}}
\put(55,65){\makebox(0,0)[cc]{$\phi$ }}
\put(83.00,54){\makebox(0,0)[cd]{ $x_1$ }}
\put(56,84){\makebox(0,0)[cd]{ $x_0$ }}
\put(50,50){\vector(1,0){35}}\put(50,50){\vector(0,1){35}}
\put(50,50){\line(-1,0){35}}\put(50,50){\line(0,-1){35}}
\put(50,75){\line(-1,0){35}}\put(50,75){\line(1,0){35}}
\put(50,25){\line(-1,0){35}}\put(50,25){\line(1,0){35}}
\end{picture}
\caption{"Rotations" or Galilei
transformations on the galileian plane}
\end{figure}
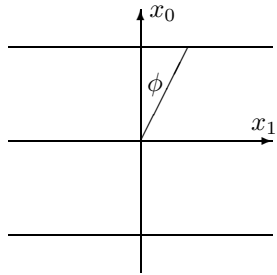

\begin{figure}[htp]
\begin{picture}(200.00,100)(-80,0)
\put(50,50){\line(1,2){11.5}}
\put(55,65){\makebox(0,0)[cc]{$\phi$ }}
\multiput(50,50)(3,3){15}{\circle*{1}}
\multiput(50,50)(-3,3){15}{\circle*{1}}
\multiput(50,50)(-3,-3){15}{\circle*{1}}
\multiput(50,50)(3,-3){15}{\circle*{1}}
\put(93.00,54){\makebox(0,0)[cd]{$ x_1$ }}
\put(56,94){\makebox(0,0)[cd]{ $x_0$ }}
\put(50,50){\vector(1,0){45}}\put(50,50){\vector(0,1){45}}
\put(50,50){\line(-1,0){45}}\put(50,50){\line(0,-1){45}}
\po(50,70,151)\po(51,70.025,152)\po(52,70.1,153)
\po(53,70.22,154)\po(54,70.4,155)\po(55,70.6,156)
\po(56,70.9,157)\po(57,71.2,158)\po(58,71.54,159)
\po(59,71.93,160)\po(60,72.36,161)\po(61,72.83,162)
\po(62,73.3,163)\po(63,73.85,164)\po(64,74.4,165)
\po(65,75,166)\po(66,75.6,167)\po(67,76.25,168)
\po(68,76.9,169)\po(69,77.6,170)\po(70,78.3,171)
\po(71,79.0,172)\po(72,79.73,173)\po(73,80.48,174)
\po(74,81.24,175)\po(75,82.02,176)\po(76,82.8,177)
\po(77,83.6,178)\po(78,84.41,179)\po(79,85.23,180)
\po(80,86.06,181)\put(80,86.06){\line(1,1){8}}
\po(49,70.025,182)\po(48,70.1,183)\po(47,70.22,184)
\po(46,70.4,185)\po(45,70.6,186)\po(44,70.9,187)
\po(43,71.2,188)\po(42,71.54,189)\po(41,71.93,190)
\po(40,72.36,191)\po(39,72.83,192)\po(38,73.3,193)
\po(37,73.85,194)\po(36,74.4,195)\po(35,75,196)
\po(34,75.6,197)\po(33,76.25,198)\po(32,76.9,199)
\po(31,77.6,200)\po(30,78.3,201)\po(29,79,202)
\po(28,79.73,203)\po(27,80.48,204)\po(26,81.24,205)
\po(25,82.02,206)\po(24,82.8,207)\po(23,83.6,208)
\po(22,84.41,209)\po(21,85.23,210)\po(20,86.06,211)
\put(20,86.06){\line(-1,1){8}}
\po(50,30,212)\po(51,29.975,213)\po(52,29.9,214)
\po(53,29.78,215)\po(54,29.6,216)\po(55,29.4,217)
\po(56,29.1,218)\po(57,28.8,219)\po(58,28.46,220)
\po(59,28.07,221)\po(60,27.64,222)\po(61,27.17,223)
\po(62,26.7,224)\po(63,26.15,225)\po(64,25.6,226)
\po(65,25,227)\po(66,24.4,228)\po(67,23.75,229)
\po(68,23.1,230)\po(69,22.4,231)\po(70,21.7,232)
\po(71,21,233)\po(72,20.27,234)\po(73,19.52,235)
\po(74,18.76,236)\po(75,17.98,237)\po(76,17.2,238)
\po(77,16.4,239)\po(78,15.6,240)\po(79,14.77,241)
\po(80,13.94,242)\put(80,13.94){\line(1,-1){8}}
\po(49,29.975,602)\po(48,29.9,603)
\po(47,29.78,604)\po(46,29.6,605)\po(45,29.4,606)
\po(44,29.1,607)\po(43,28.8,608)\po(42,28.46,609)
\po(41,28.07,610)\po(40,27.74,611)\po(39,27.17,612)
\po(38,26.7,613)\po(37,26.15,614)\po(36,25.6,615)
\po(35,25,616)\po(34,24.4,617)\po(33,23.75,618)
\po(32,23.1,619)\po(31,22.4,620)\po(30,21.7,621)
\po(29,21,622)\po(28,20.27,623)\po(27,19.52,624)
\po(26,18.76,625)\po(25,18,626)\po(24,17.2,627)
\po(23,16.4,628)\po(22,15.6,629)\po(21,14.77,630)
\put(21,14.77){\line(-1,-1){8}}
\con(151,152)\con(152,153)\con(153,154)\con(154,155)
\con(155,156)\con(156,157)\con(157,158)\con(158,159)
\con(159,160)\con(160,161)\con(161,162)\con(162,163)
\con(163,164)\con(164,165)\con(165,166)\con(166,167)
\con(167,168)\con(168,169)\con(169,170)\con(170,171)
\con(171,172)\con(172,173)\con(173,174)\con(174,175)
\con(175,176)\con(176,177)\con(177,178)\con(178,179)
\con(179,180)\con(180,181)\con(151,182)\con(182,183)
\con(183,184)\con(184,185)\con(185,186)
\con(186,187)\con(187,188)\con(188,189)
\con(189,190)\con(190,191)\con(191,192)
\con(192,193)\con(193,194)\con(194,195)
\con(195,196)\con(196,197)\con(197,198)
\con(198,199)\con(199,200)\con(200,201)
\con(201,202)\con(202,203)\con(203,204)
\con(204,205)\con(205,206)\con(206,207)
\con(207,208)\con(208,209)\con(209,210)
\con(210,211)\con(212,213)\con(213,214)
\con(214,215)\con(215,216)\con(216,217)
\con(217,218)\con(218,219)\con(219,220)
\con(220,221)
\con(221,222)\con(222,223)\con(223,224)\con(224,225)
\con(225,226)\con(226,227)\con(227,228)\con(228,229)
\con(229,230)\con(230,231)
\con(231,232)\con(232,233)\con(233,234)
\con(234,235)\con(235,236)\con(236,237)
\con(237,238)\con(238,239)\con(239,240)
\con(240,241)\con(241,242)
\con(212,602)\con(602,603)
\con(603,604)\con(604,605)\con(605,606)
\con(606,607)\con(607,608)\con(608,609)\con(609,610)
\con(610,611)\con(611,612)\con(612,613)\con(613,614)
\con(614,615)\con(615,616)\con(616,617)\con(617,618)
\con(618,619)\con(619,620)\con(620,621)\con(621,622)
\con(622,623)\con(623,624)\con(624,625)\con(625,626)
\con(626,627)\con(627,628)\con(628,629)\con(629,630)
\end{picture}
\caption{Hyperbolic rotations or Lorentz transformations on the minkowskian plane}
\end{figure}
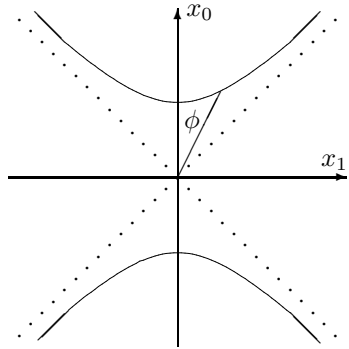

 Two dimentional vectors
$ x^t=(x_0,x_1)^t $ in Cartesian basis are transformed
under rotations as follows:
\begin{equation}
x'=Ax,
\label{2.1k}
\end{equation}
where
$$
\displaystyle{A=\left ( \begin{array}{cc}
\cos\varphi & -\sin\varphi \\
\sin\varphi & \cos\varphi \end{array} \right ) \in SO(2)}
$$
is the ordinary rotation matrix,
$$
 \displaystyle{A=\left ( \begin{array}{cc}
1 & 0 \\
\varphi & 1 \end{array} \right ) \in E(1)}
$$
is the Galilei transformation and
$$
 \displaystyle{A=\left ( \begin{array}{cc}
\cosh\varphi & \sinh\varphi \\
\sinh\varphi & \cosh\varphi \end{array} \right ) \in SO(1,1)}
$$
is the hyperbolic rotation or Lorentz transformation, respectively.
The following quadratic forms are invariant under
transformations (\ref{2.1k}): $ inv_e=x_0^2+x_1^2, \ inv_g=x_0^2,
 \ inv_m=x_0^2-x_1^2. $
The Lie algebras of these three groups may be written in a unified manner
\begin{equation}
 \left ( \begin{array}{cc} o & \omega \\
1 & 0 \end{array} \right ) \in
so(2;\omega),
\ \ \omega=1, 0, -1,
\label{2.2k}
\end{equation}
where $ \omega=1 $ correspond to Lie algebra of the simple rotation group
$ SO(2), \omega=-1 $ correspond to the semisimple
Lorentz group $ SO(1,1) $ and $ \omega=0 $ correspond to the
nonsemisimple Galilei group $ E(1). $

The one dimensional geometries of constant curvature
are realized on the rotation invariant surfaces
(spheres) in these planes. Eccording with Erlanger program due to
F.Klein a geometry is completely
determined by its motion group. In one dimensional case there is only
one motion, namely translation. Let us introduce the intrinsic
(Beltrami) coordinate on the spheres by relation
$ \xi=x_1/x_0. $ Then the translation operators
$ T(a): \xi \rightarrow \xi'=x'_1/x'_0 $ looks as follows:
\begin{equation}
\xi'=T(a)\xi=\frac{\xi+a}{1-a\xi}, \quad  a=\tan\varphi,
\quad a \in {\bf R}
\label{2.3}
\end{equation}
for the elliptic (spherical) geometry of constant positive curvature,
\begin{equation}
\xi'=T(a)\xi=\xi+a, \quad a=\varphi, \quad a \in {\bf R}
\label{2.4}
\end{equation}
for the euclidean flat geometry and
\begin{equation}
\xi'=T(a)\xi=\frac{\xi+a}{1+a\xi}, \quad  a=\tanh\varphi,
\quad a \in (-1,1)
\label{2.5}
\end{equation}
for the hyperbolic geometry of constant negative curvature.
Translation operators may be written in the following unified manner
\begin{equation}
\xi'=T(a;\omega)\xi=\frac{\xi+a}{1-\omega{a\xi}}, \quad
 a=\frac{1}{\sqrt{\omega}}\tan\sqrt{\omega}\varphi,
\label{2.6}
\end{equation}
with the same parameter $ \omega $ as in Eq.(\ref{2.2k}).
The distance $ d_{AB} $ between points $ A $ and $ B $
\begin{equation}
\frac{1}{\sqrt{\omega}}\tan(\sqrt{\omega}d_{AB})=
\left| \frac{\xi_B-\xi_A}{1+\omega\xi_B\xi_A} \right|
\label{2.7}
\end{equation}
is invariant under translations.

The new parameter $ j=\sqrt{\omega} $ is appeared in
 Eqs.(\ref{2.6}),(\ref{2.7}),
where $ j=1 $ for $ \omega=1, \ j=i $ for $ \omega=-1 $
and $ j=\sqrt{0} $ for $ \omega=0. $
The solutions of equation $ j=\sqrt{0} $ depends on underlying mathematical
structure, for example $ j=0 $ over fields {\bf R} or {\bf C}; \
$ \displaystyle{ j=\left ( \begin{array}{cc}
0 & 0 \\
1 & 0 \end{array} \right ) \neq 0} $ over $ 2\times2 $ matrix ring;
$ j_k=\theta_k, \  \theta_k^2=0, \  \theta_k\theta_m=-\theta_m\theta_k,
\ k \neq m, \ k,m=1, \ldots, N $ over Grassmann algebra.
We take as the solution nilpotent commutative numbers
$ j_k=\iota_k, \ \iota_k^2=0, \  \iota_k \iota_m=\iota_m \iota_k, \
k \neq m, \ k,m=1,\ldots,N. $

For $ k=1 $ the nilpotent number $ \iota, \  \iota^2=0, $
have been first introduced by W.K.Clifford [4] more then
hundred years ago and was applied in geometry and mechanics
[5]---[7]. In Russian publications its is named as dual number and in
English is known as Study number. R.I.Pimenov [8] was the first who have
introduce the {\it different} nilpotent numbers $ \iota_k, \ k=1, \ldots,N $
with {\it commutative} law of multiplication, so it seems natural
to call its as Pimenov numbers.

Now we are able to rewrite transformations
$ x'=A(j)x $ in the form
\begin{equation}
\left ( \begin{array}{c}
x'_0 \\
x'_1 \end{array} \right )=
\left ( \begin{array}{cc}
\cos j\varphi & -j\sin j\varphi \\
j^{-1}\sin j\varphi & \cos j\varphi \end{array} \right )
\left ( \begin{array}{c}
x_0 \\
x_1 \end{array} \right ),
\label{2.8k}
\end{equation}
where $ \det{A(j)}=\cos^2j\varphi+\sin^2j\varphi=1, \
inv(j)=x^2_0+j^2x_1^2 $ and parameter $ j $ takes three values
$ j=1, \iota, i. $
Matrices $ A(j) $ form the group $ SO(2;j). $
It is easily to check that $ SO(2;1)\equiv SO(2), \
SO(2;i)\equiv SO(1,1) $ and
$ SO(2;\iota) \equiv E(1). $ In the last case the
properties of a functions of $ \iota $ arising from its
Taylor expansion are exploited:
$ \cos\iota\varphi=1, \  \sin\iota\varphi=\iota\varphi $
and be definition $ \iota/\iota=1. $
If we put $ j=\epsilon \in {\bf R}, \  \epsilon \rightarrow 0, $
then we obtain In\"on\"u-Wigner contraction [9] on the group level.
The remarkable property of Eq.(\ref{2.8k}) is that the matrix elements
and vector components are {\it real} numbers for any
value of the parameter $ j.$

There is other way to describe rotations (\ref{2.1k}),
namely $ x'(j)=R(j)x(j) $ or more precisely
\begin{equation}
\left ( \begin{array}{c}
x'_0 \\
jx'_1 \end{array} \right )=
\left ( \begin{array}{cc}
\cos j\varphi & -\sin j\varphi \\
\sin j\varphi & \cos j\varphi \end{array} \right )
\left ( \begin{array}{c}
x_0 \\
jx_1 \end{array} \right ),
\label{2.9}
\end{equation}
where $ \det{R(j)}=\cos^2j\varphi+\sin^2j\varphi=1 $ and
$ inv(j)=x^2_0+j^2x_1^2. $
Matrices $ R(j) $ again form the same group $ SO(2;j). $
The remarkable property of Eq.(\ref{2.9}) is that some matrix
elements and vector components are
{\it nilpotent Pimenov} numbers for $ j=\iota $
\begin{equation}
\left ( \begin{array}{c}
x'_0 \\
{\iota}x'_1 \end{array} \right )=
\left ( \begin{array}{cc}
1  & -{\iota}\varphi \\
{\iota}\varphi & 1 \end{array} \right )
\left ( \begin{array}{c}
x_0 \\
{\iota}x_1 \end{array} \right ),
\label{2.10}
\end{equation}
or {\it complex} numbers for $ j=i $
\begin{equation}
\left ( \begin{array}{c}
x'_0 \\
ix'_1 \end{array} \right )=
\left ( \begin{array}{cc}
\cosh\varphi & -i\sinh\varphi \\
i\sinh\varphi & \cosh\varphi \end{array} \right )
\left ( \begin{array}{c}
x_0 \\
ix_1 \end{array} \right ).
\label{2.11}
\end{equation}

This simple considiration suggests to regard Lie groups,
Lie algebras and vector spaces not only over complex
(or real) number field, but also over more general algebraical
structure (see next Section) which is generated by
nilpotent commutative numbers.

Let us observe that matrix
$ R(j) $ in Eq.(\ref{2.9}) is obtained from rotation matrix
$ A \in SO(2) $ by substitution of the real group parameter
$ \varphi $ with new group parameter
$ j\varphi. $ For $ j=i $ this
means the analytical continuation of the real group
parameter into the complex number field
(Weyl's unitary trick), which gives in result the noncompact
group $ SO(1,1) $ from the simple group $ SO(2). $
In the same line In\"on\"u-Wigner contraction
$ (j=\iota) $ may be regarded as the continuation of the real
$ \varphi $ to the nilpotent values $ \iota\varphi. $

In the traditional approach
$ j=\epsilon, $ Eq.(\ref{2.9}) looks as follows
\begin{eqnarray}
x'_0 & = & x_0\cos\epsilon\varphi-\epsilon x_1\sin\epsilon\varphi
\nonumber \\
\epsilon x'_1 & =
& \epsilon x_1\cos\epsilon\varphi+x_0\sin\epsilon\varphi
\label{2.12}
\end{eqnarray}
and in the limit $ \epsilon \rightarrow 0 $ provide
the Galilei transformation if the infinitesimals of the first order
are compared in both sides of the second equation. Let us stress
that the group transformation
$ R(\epsilon) $ and vector $ x(\epsilon) $ are regarded
together in Eq.(\ref{2.12}) and can not be separated.
Otherwise for example
$ \lim_{\epsilon \rightarrow 0}R(\epsilon)=I $
instead of to be a general element of the Galilei group
$ E(1). $ Just the similar situation is appeared in the case of
quantum groups (see for example [13]) and it seems
that the  nilpotent
commutative numbers provide
the relevant mathematical structure to avoid such troubles.

\section{ Pimenov algebra $ {\bf D}_{n}(\iota;{\bf C}) $ }

     Algebra $ {\bf D}_{n}(\iota;{\bf C}) $ is defined as an associative
  algebra with unit and {\it nilpotent} generators
   $ \iota_{1},\ldots,\iota_{n},\; \iota_{k}^{2}=0,\; k=1,\ldots,n $
 with {\it commutative} multiplication
     $ \iota_{k}\iota_{m}=\iota_{m}\iota_{k},\; k\neq{m}. $
  The general element of $ {\bf D}_{n}(\iota;{\bf C}) $ has the form
\begin{equation}
a=a_0 + \sum_{p=1}^{n} \sum_{k_1<...<k_p}
a_{k_1...k_p}\iota_{k_1}...\iota_{k_p}, \quad
a_0,a_{k_1...k_p}\in{\bf C}.
\end{equation}
For $ n=1 $ we have
$ {\bf D}_{1}(\iota_{1};{\bf C})\ni{a=a_{0}+a_{1}\iota_{1}}, $
i.e. dual (or Study) numbers, when
$ a_{0},a_{1}\in{\bf R}. $
For $ n=2 $ the general element of
$ {\bf D}_{2}(\iota_{1},\iota_{2};{\bf C}) $
is written as follows:
$ a=a_{0}+a_{1}\iota_{1}+a_{2}\iota_{2}+a_{12}\iota_{1}\iota_{2}. $

    Two elements $ a,b \in   {\bf D}_{n}(\iota;{\bf C}) $ are
equal    when
 $  a_0=b_0, \  a_{k_1...k_p}=b_{k_1...k_p}, \  p=1,\ldots,n. $
If $ a=a_k\iota_k $  and $ b=b_k\iota_k, $ then the condition
$ a=b, $ which is equivalent to $ a_k\iota_k=b_k\iota_k $
make possible to define consistently divivision of Pimenov
unit $ \iota_k $ by itself, namely $ \iota_k / \iota_k =1. $
 Divisions of a real or complex numbers by Pimenov units
 $ z/{\iota_k}, \  z \in {\bf R}, {\bf C},$ and different
  Pimenov units
$ {\iota_m}/{\iota_k}, \  k \neq m $ are not defined.
 A function
 $ f: {\bf D}_{n}(\iota;{\bf C}) \rightarrow
 {\bf D}_{n}(\iota;{\bf C}) $
 is defined by its Taylor expansion
      $$
f(a) = f(a_0) + \sum_{p=1}^{n} \sum_{k_1<...<k_p}
f_{k_1...k_p}\iota_{k_1}...\iota_{k_p},
      $$
\begin{equation}
f_{k_1...k_p}=\sum_{r=1}^{n} f^{(r)}(a_0)d(p;r), \quad
d(p;r)=\sum_{v_1 \cup \ldots \cup v_r =(k_1, \ldots ,k_p)}
      a_{v_1} \cdots a_{v_r},
\label{3.2}
\end{equation}
 where  the summation on all possible
 partitions of the number set $ (k_1, \ldots ,k_p) $ on $ p $
 nonempty subsets is understood in the last equation.
    For example,
    $$
    a\in{\bf D}_{1}(\iota;{\bf C}), \quad
 f(a) = f(a_0)+{\iota_1}a_1{f^{'}(a_0)};
 $$
      $$
    a\in{\bf D}_{2}(\iota;{\bf C}), \quad
 f(a)=f(a_0)+{\iota_1}a_1{f^{'}(a_0)}+{\iota_2}a_2{f^{'}(a_0)}
 +{\iota_1}{\iota_2}(a_{12}{f^{'}(a_0)+a_1a_2f^{''}(a_0))}.
 $$
      It is clear from previous equations that a function
 over algebra
 $ {\bf D}_{n}(\iota;{\bf C}) $
   is completely determined by its real part $ f(a_0) $
   and $ n $ derivatives $ f^{(r)}(a_0), \; r=1,\ldots,n, $
 for example over
 $ {\bf D}_{1}(\iota;{\bf C}) $
$$
      e^{a_0+\iota_{1}a_1} = e^{a_0} + \iota_{1}a_{1}e^{a_0} =
      e^{a_0}(1 + \iota_{1}a_{1}),
$$
\begin{equation}
      \sin (a_0 + \iota_{1}a_{1}) =
      \sin a_0 + \iota_{1}a_{1}\cos a_0,
\label{3.3}
\end{equation}
 and over
 $ {\bf D}_{2}(\iota;{\bf C}) $
$$
      e^{a} = e^{a_0}(1 + \iota_{1}a_{1} + \iota_{2}a_{2}
 + \iota_{1}\iota_{2}(a_{12}+a_1a_2)),
 $$
 \begin{equation}
      \sin{(a)} = \sin{a_0} + (\iota_{1}a_{1}
                + \iota_{2}a_{2})\cos{a_0}
  +  \iota_{1}\iota_{2}(a_{12}\cos{a_0} - a_1a_2\sin{a_0}).
\label{3.4}
\end{equation}

The well known Grassmann algebra
   $ {\bf \Gamma}_{n}( \xi ) $
 is the algebra with {\it nilpotent} generators
     $ \xi_k^2=0, k=1,\ldots,n $
 and {\it anticommutative} multiplication
     $ \xi_k\xi_m= -\xi_m\xi_k, k \neq m. $
  It is easy to verify that the product of two generators of
 Grassmann algebra has the same algebraic properties as the
 generator of  algebra $ {\bf D}_{n}(\iota;{\bf C}): $
     $ \iota_k=\xi_k\xi_{n+k}, \; k=1,\ldots,n. $
  This means that this  algebra is the subalgebra of even
  part of Grassmann algebra
     $ {\bf D}_{n}(\iota;{\bf C}) \subset   {\bf \Gamma}_{2n}(\xi ). $

For our aims it is convenient to regard the set of algebras
$ {\bf D}_n(j;{\bf C}) $
with $ j_k=1,\iota_k, \; k=1,\ldots,n. $
If some parameters are equil to Pimenov numbers
$ j_{k_s}=\iota_s, \; s=1,\ldots,m $ and remaining
ones are equal to 1,  then we have the algebra
$ {\bf D}_m(\iota;{\bf C}) $
from the set
$ {\bf D}_n(j;{\bf C}). $

\section{ Orthogonal CK groups  }

Let us regard according to R.I.Pimenov [8] a specific
vector space ${\bf R}_{N}(j) $ over
$ {\bf D}_{N-1}(j;{\bf R}) $
with Cartesian coordinates
$ x(j)=(x_1,J_{12}x_2,\ldots, J_{1,N}x_{N})^t, \; x_k\in{\bf R},
 \;k=1,\ldots,N $
and quadratic form
\begin{equation}
x^t(j)x(j) = x_1^2 + \sum_{k=2}^{N}J_{1k}^2x_k^2,
\label{2.1}
\end{equation}
where
     \begin{equation}
 J_{\mu\nu}=\prod_{r=\mu}^{\nu -1}j_{r}, \quad
     \mu < \nu, \quad  J_{\mu\nu}=1, \quad \mu \geq \nu,
     \quad  j_{r}=1,\iota_{r},i.
    \label{2.2}
    \end{equation}
 Orthogonal CK groups $ SO(N;j;{\bf R}) $
 are defined as the
 set of transformations of $ {\bf R}_{N}(j) $ leaving invariant
 (\ref{2.1}) and are
 realized in the
 Cartesian basis as the matrix groups over ${\bf D}_{N-1}(j;{\bf R}) $
 with the help of the {\it special} matrices
     \begin{displaymath}
   (A(j))_{kp}=\tilde{J}_{kp}a_{kp}, \quad a_{kp}\in{\bf R},
    \end{displaymath}
     \begin{equation}
  \tilde{J}_{kp}=J_{kp},\; k<p,   \quad    \tilde{J}_{kp}=J_{pk},\;
k\geq{p}.
    \label{16}
    \end{equation}
 These matrices  act on  vectors  $ x(j)\in{\bf R}_{N}(j) $
 by matrix multiplication and
  are   satisfied   the  following $ j$-orthogonality
relations:
     \begin{equation}
    A(j)A^{t}(j)=A^{t}(j)A(j)=I.
    \label{17}
    \end{equation}

   Let in euclidean vector space $ {\bf R}_n \; y=Dx $
   is the transformation from Cartesian basis $ x $ to the
   new ("symplectic") basis $ y $ with the following
   quadratic form: $ y^tC_0y, $ where
   $ (C_0)_{ik}=\delta_{ik'},\; k'=N+1-k. $
   The matrix $ D $ is obtained from the invariant condition
   for quadratic form, i.e.
   $ y^tC_0y=x^tD^tC_0Dx=x^tx $ for any $ x \in {\bf R}_n. $
     One of the solutions of matrix equation
    $ D^tC_0D=I $
     has the form
\begin{equation}
D=\frac{1}{\sqrt{2}}
\left ( \begin{array}{cc}
      I &      {\tilde C_0} \\
      i{\tilde C_0} &   -iI
      \end{array} \right ),    \
D=\frac{1}{\sqrt{2}}
\left ( \begin{array}{ccc}
      I & 0 &   {\tilde C_0} \\
      0 & \sqrt{2} &  0 \\
      i{\tilde C_0} & 0  &  -iI
      \end{array} \right ),    \
\label{8}
\end{equation}
for $ N=2n $ and $ N=2n+1, $ respectively,
    where $ \tilde{C}_{0}\in{M_{n}({\bf C})} $ is the matrix with
 the real units on the second diagonal.

The similarity transformation
\begin{equation}
 B(j)=D^{-1}A(j)D
\label{18}
\end{equation}
in the vector space $ {\bf R}_n(j)  $
gives the realization of
$ SO(N;j;{\bf R}) $  in a new ("symplectic") basis
$ y(j)=Dx(j) $
with the invariant quadratic form
\begin{equation}
y^t(j)C_0y(j) = 2\sum_{k=1}^{n}J_{1k}J_{1k'}y_{k}y_{k'} +
\epsilon J^2_{1,n+1}y^2_{n+1},
\label{2.8}
\end{equation}
$ \epsilon=1 $ for $ N=2n+1, \; \epsilon=0 $ for $ N=2n $
and  the additional relations of $ j$-orthogonality
     \begin{equation}
   B(j)C_{0}B^{t}(j)=B^{t}(j)C_{0}B(j)=C_{0}.
    \label{19}
    \end{equation}

The solution $ D $ (\ref{8}) is not unique.
There are different solutions $ \tilde{D} $ of the matrix
equation $ D^tC_0D=I. $
The similarity transformations
$ \tilde{B}(j)=\tilde{D}^{-1}A(j)\tilde{D} $
provide the different realization of
$ SO(N;j;{\bf R}) $ as the matrix group over
$ {\bf D}_{N-1}(j;{\bf R}).$
In the case of quantum groups
$ SO_q(N;j;{\bf C}) $
this correspond to the different couplings of CK and Hopf structures
on the level of quantum groups, which on the
level of quantum algebras mean the different choice of the
primitive elements of the Hopf algebra [15].

\section{ Quantum orthogonal CK groups $ SO_{v}(N;j;{\bf C}) $ }

     According to FRT theory  of quantum groups  the starting point
     of quantization is
 an algebra $ {\bf C} \langle t_{ik} \rangle $
of noncommutative polynomials of
$ N^2 $ variables $ t_{ik}, i,k=1,\ldots,N $ over complex number
field ${\bf C}. $ For well known [1] lower triangular matrix
$ R_{q}\in{M_{N^2}({\bf C})} $ the generators
$ T=(t_{ik})_{i,k=1}^{N}\in{M_{N}({\bf C} \langle t_{ik} \rangle)} $
have the following commutation relations
\begin{equation}
R_qT_1T_2=T_2T_1R_q,
\label{2}
\end{equation}
where $ T_1=T \otimes I, \  T_2=I \otimes T \in M_{N^2}({\bf C}
\langle t_{ij} \rangle). $
There are additional relations of $ q$-orthogonality
\begin{equation}
TCT^t=T^tCT=C,
\label{3}
\end{equation}
where $ C=C_0q^{\rho},  \rho=diag(\rho_1, \ldots, \rho_N), $
\begin{equation}
(\rho_1, \ldots, \rho_N)=
\left \{ \begin{array}{cc}
     (n-\frac{1}{2}, n-\frac{3}{2}, \ldots , \frac{1}{2},0,-\frac{1}{2},
     \ldots , -n+\frac{1}{2}), &  N=2n+1 \\
     (n-1, n-2, \ldots, 1,0,0,-1, \ldots, -n+1), &  N=2n.
     \end{array} \right.
\label{4}
\end{equation}
The quantum orthogonal group $ SO_{q}(N;{\bf C}) $
is defined as the quotient
     \begin{equation}
SO_{q}(N;{\bf C})={\bf C} \langle t_{ik} \rangle \big/ (\ref{2}),(\ref{3}).
    \label{5}
    \end{equation}
>From the algebraic point of view $ SO_{q}(N;{\bf C}) $ is
a Hopf algebra with the following coproduct $ \Delta, $
counit $ \epsilon $ and antipode $ S: $
\begin{equation}
\Delta T=T \dot {\otimes} T,\quad \epsilon(T)=I,\quad S(T)=CT^tC^{-1}.
\label{6}
\end{equation}

     We shall regard the quantum deformations of the contracted complex CK
 groups and in this case the parameters $ j $ take only two values:
  $ j_{k}=1,\iota_{k}. $
 We start now with the
   $ {\bf D}\langle t_{ik} \rangle  $  ---  the algebra of  noncommutative
 polynomials of $ N^2 $ variables  over
 the  algebra $ {\bf D}_{N-1}(j). $ In addition we  transform
 the deformation parameter $ q=\exp{z} $ as follows:
     \begin{equation}
     z=Jv, \quad J\equiv J_{1N}=\prod_{k=1}^{N-1}j_{k},
    \label{20}
    \end{equation}
 where $ v $ is the new deformation parameter. The
 transformation of quantum deformation parameter was suggested
 by E.Celeghini et al. [11].

      In "symplectic" basis the quantum CK group $ SO_{v}(N;j;{\bf C}) $
 is produced by the generating matrix
$ T(j)\in{M_{N}({\bf D} \langle t_{ik} \rangle ) } $
 equal to $ B(j) $ (\ref{18}) for $ q=1. $
 The non\-com\-mu\-ta\-ti\-ve entries of
 $ T(j) $ obey the commutation relations
\begin{equation}
 R_v(j)T_1(j)T_2(j)=T_2(j)T_1(j)R_v(j)
\label{21}
\end{equation}
 and the additional relations of $ (v,j)$-orthogonality
     \begin{equation}
   T(j)C(j)T^{t}(j)=T^{t}(j)C(j)T(j)=C(j),
    \label{22}
    \end{equation}
 where lower triangular R--matrix $ R_{v}(j) $ and $ C(j) $ are
 obtained from $ R_q $ and $ C, $ respectively, by substitution
 $ Jv $ instead of $ z: $
     \begin{equation}
     R_{v}(j)=R_{q}(z \rightarrow Jv), \quad
     C(j)=C(z \rightarrow Jv).
    \label{23}
    \end{equation}
 Then the quotient
     \begin{equation}
SO_{v}(N;j;{\bf C})=
 {\bf D} \langle t_{ik} \rangle \big/ (\ref{21}),(\ref{22})
    \label{24}
    \end{equation}
is Hopf algebra with the  coproduct $ \Delta, $
counit $ \epsilon $ and antipode $ S: $
\begin{equation}
\Delta T(j)=T(j) \dot {\otimes}T(j),\quad \epsilon (T(j))=I, \quad
 S(T(j))=C(j)T^t(j)C^{-1}(j).
\label{25}
\end{equation}

The matrix $ D $ (\ref{8}) in the similarity transformation
(\ref{18}) and the factor
$ J=J_{1N} $ in the deformation parameter
transformation (\ref{20}) are selected consistently to provide
the existence of the Hopf algebra structure
for all possible values of the parameters
$ j, $ i.e. for all contracted CK groups. For some other
solution $ \tilde{D} $ in Eq.(\ref{18}) the consistent
factor $ J $ in Eq.(\ref{20}) may be equal to product only some
parameters $ j_k. $ It turn out that for some choice of
$ \tilde{D} $ not all CK contractions are allowed.

\section{ Quantum CK algebras $ so_{v}(N;j;{\bf C}) $ as a dual to
     $ SO_{v}(N;j;{\bf C}) $ }

 By FRT quantization  theory [1] the dual space
   $ Hom(SO_{v}(N;j;{\bf C}),{\bf C}) $ is an algebra with the multiplication
 induced by coproduct $ \Delta $ in $ SO_{v}(N;j;{\bf C}) $
     \begin{equation}
     l_{1}l_{2}(a)=(l_{1}\otimes l_{2})(\Delta (a)),
    \label{30}
    \end{equation}
  $ l_{1},l_{2}\in{ Hom(SO_{v}(N;j;{\bf C}),{\bf C})},
 \quad a\in{SO_{v}(N;j)}. $
 Let us formally introduce $ N \times N $ upper $ (+) $ and lower
 $ (-) $ triangular matrices $ L^{(\pm)}(j) $ as follows: it is
 necessary to put $ j_{k}^{-1} $ in the nondiagonal matrix elements
  of $ L^{(\pm)}(j), $ if there is the parameter $ j_k $ in the
  corresponding matrix element of $ T(j). $ For example, if
  $ (T(j))_{12}=j_{1}t_{12}+j_{2}\tilde{t}_{12}, $ then
  $ (L^{(+)}(j))_{12}=j_{1}^{-1}l_{12}+j_{2}^{-1}\tilde{l}_{12}. $
  Formally the matrices   $   L^{(\pm)}(j)   $   are   not   defined  for
  $ j_{k}=\iota_{k}, $ since $ \iota_{k}^{-1} $ do not exist, but
  $ l_{ik}^{(\pm)} $ are functionals on $ t_{pr}, $ so
  if we set an action of the matrix functionals $ L^{(\pm)}(j) $
  on the elements of $ SO_{v}(N;j;{\bf C}) $ by the duality relation
     \begin{equation}
    \langle L^{(\pm)}(j),T(j) \rangle = R^{(\pm)}(j),
    \label{31}
    \end{equation}
where
     \begin{equation}
     R^{(+)}(j)=PR_{v}(j)P, \quad  R^{(-)}(j)= R_{v}^{-1}(j),\quad
     Pu \otimes w = w \otimes u,
    \label{32}
    \end{equation}
then we shall have well defined expressions even for
$ j_{k}=\iota_k. $

     The elements of $ L^{(\pm)}(j) $ satisfy the commutation
 relations
     \begin{eqnarray}
     R^{(+)}(j)L_{1}^{(\sigma)}(j)L_{2}^{(\sigma)}(j) & = &
 L_{2}^{(\sigma)}(j)L_{1}^{(\sigma)}(j)R^{(+)}(j), \nonumber  \\
     R^{(+)}(j)L_{1}^{(+)}(j)L_{2}^{(-)}(j) & = &
L_{2}^{(-)}(j)L_{1}^{(+)}(j)R^{(+)}(j),  \quad
     \sigma = \pm
    \label{33}
    \end{eqnarray}
and additional relations
     \begin{eqnarray}
  L^{(\pm)}(j)C^{t}(j)(L^{(\pm)}(j))^{t} & = & C^{t}(j),  \nonumber \\
(L^{(\pm)}(j))^{t}(C^{t}(j))^{-1}L^{(\pm)}(j)  & =
 & (C^{t}(j))^{-1}, \nonumber \\
l_{kk}^{(+)}l_{kk}^{(-)}=l_{kk}^{(-)}l_{kk}^{(+)}=1, & &
     l_{11}^{(+)}\ldots l_{NN}^{(+)}=1, \; k=1,\ldots ,N.
    \label{34}
    \end{eqnarray}
 An algebra  $  so_{v}(N;j;{\bf C})=\{I,L^{(\pm)}(j)\} $ is called quantum CK
algebra and
is Hopf algebra with the following  coproduct $ \Delta, $
counit $ \epsilon $ and antipode $ S: $
  $$
 \Delta L^{(\pm)}(j)=L^{(\pm)}(j) \dot {\otimes} L^{(\pm)}(j), \quad
\epsilon (L^{(\pm)}(j))=I,
  $$
\begin{equation}
S(L^{(\pm)}(j))  =  C^{t}(j)(L^{(\pm)}(j))^{t}(C^{t}(j))^{-1}.
    \label{35}
    \end{equation}

 It is possible to show that algebra $ so_{v}(N;j;{\bf C}) $ is isomorphic
 with the quantum deformation [10] of the universal enveloping algebra
 of the CK algebra $ so(N;j;{\bf C}), $ which may be obtained from the
 orthogonal algebra $ so(N;{\bf C}) $ by contractions [2].
So there are at least two ways for construction of quantum CK
algebras.

\section{ Example: $ SO_{v}(3;j;{\bf C}) $ and $ so_{v}(3;j;{\bf C}) $ }

The generating matrix for the simplest quantum orthogonal group
  $ SO_{v}(3;j;{\bf C}),  \; j=(j_1,j_2) $ is in the form
     \begin{equation}
        T(j) =
     \left( \begin{array}{ccc}
 t_{11}+ij_{1}j_{2}\tilde{t}_{11} & j_{1}t_{12}-ij_{2}\tilde{t}_{12}
& t_{13}-ij_{1}j_{2}\tilde{t}_{13} \\
j_{1}t_{21}+ij_{2}\tilde{t}_{21}  & t_{22}
& j_{1}t_{21}-ij_{2}\tilde{t}_{21}       \\
 t_{13}+ij_{1}j_{2}\tilde{t}_{13} & j_{1}t_{12}+ij_{2}\tilde{t}_{12}
& t_{11}-ij_{1}j_{2}\tilde{t}_{11} \\
     \end{array} \right) .
    \label{36}
    \end{equation}
The R-matrix is obtained from the standart one by Eq.(\ref{23})
and is as follows
$$
        R_{v}(j)\equiv R_{q}(z \rightarrow Jv) =
$$
     \begin{equation}
     \left( \begin{array}{ccccccccc}
e^{Jv}  & 0 & 0 & 0 & 0 & 0 & 0 & 0 & 0  \\
     0  & 1 & 0 & 0 & 0 & 0 & 0 & 0 & 0  \\
     0  & 0 & e^{-Jv} & 0 & 0 & 0 & 0 & 0 & 0  \\
     0  & 2\sinh Jv  & 0 & 1 & 0 & 0 & 0 & 0 & 0  \\
          0  & 0 & -2e^{-Jv/2}\sinh Jv & 0 & 1 & 0 & 0 & 0 & 0  \\
          0  & 0 & 0 & 0 & 0 & 1 & 0 & 0 & 0  \\
          0  & 0 & 2(1-e^{-Jv})\sinh Jv & 0 & -2e^{-Jv/2}\sinh Jv & 0
 & e^{-Jv}  & 0 & 0  \\
          0  & 0 & 0 & 0 & 0 & 2\sinh Jv & 0 & 1 & 0  \\
          0  & 0 & 0 & 0 & 0 & 0 & 0 & 0 & e^{Jv}  \\
     \end{array} \right) ,
    \label{39}
    \end{equation}
   where  $ J=j_1j_2. $
 Over the  algebras $ {\bf D}_2(j_1,j_2),\; j_1=\iota_1, j_2=1, $ or $
j_1=1, j_2=\iota_2, $ or $  j_1=\iota_1, j_2=\iota_2 $ this R-matrix
  may be written in the form
     \begin{equation}
   R_{v}(j) = I + Jv \tilde R,
    \label{40}
    \end{equation}
where
     \begin{equation}
     (\tilde R)_{11}=(\tilde R)_{99}=1, \
     (\tilde R)_{33}=(\tilde R)_{77}=-1, \
     (\tilde R)_{42}=(\tilde R)_{86}=2, \
     (\tilde R)_{53}=(\tilde R)_{75}=-2
    \label{41}
    \end{equation}
    and all other elements of the matrix $ \tilde R $
    are equal to zero.
The commutation  relations  and  additional  relations  of  $ (v,j)
$-orthogonality may be obtained from Eqs.(\ref{21}),(\ref{22}) by
straitforward calculations, so we shall concentrate our attention
on the construction of quantum algebra $ so_{v}(3;j;{\bf C}). $

     The matrix functionals $ L^{\pm}(j) $ have the form
     \begin{equation}
        L^{(+)}(j) =
     \left( \begin{array}{ccc}
 l_{11} & j_{1}^{-1}l_{12}-ij_{2}^{-1}\tilde{l}_{12}
& l_{13}-ij_{1}^{-1}j_{2}^{-1}\tilde{l}_{13} \\
0  & 1 & j_{1}^{-1}l_{21}-ij_{2}^{-1}\tilde{l}_{21}       \\
 0 & 0 & l_{11}^{-1} \\
     \end{array} \right) ,
    \label{37}
    \end{equation}
     \begin{equation}
        L^{(-)}(j) =
     \left( \begin{array}{ccc}
l_{11}^{-1}  & 0 & 0 \\
-(j_{1}^{-1}l_{21}+ij_{2}^{-1}\tilde{l}_{21})  & 1 & 0       \\
-(l_{13}+ij_{1}^{-1}j_{2}^{-1}\tilde{l}_{13})  &
-(j_{1}^{-1}l_{12}+ij_{2}^{-1}\tilde{l}_{12})
&  l_{11}  \\
     \end{array} \right) .
    \label{38}
    \end{equation}
Their actions on the generators (\ref{36}) of quantum group
  $ SO_{v}(3;j;{\bf C}) $ are given by Eq.(\ref{31}) and are as follows [3]:
     $$
  l_{11}(t_{22})=1, \quad l_{11}(t_{11})=\cosh Jv, \quad
  l_{11}(\tilde{t}_{11})=-J^{-1}\sinh Jv, \quad
     $$
     $$
  l_{12}(\tilde{t}_{21})=-ij_1^2J^{-1}\sinh Jv,           \quad
  l_{12}(\tilde{t}_{12})=ij_1^2(2J)^{-1}(\sinh 3Jv/2 + \sinh Jv/2 ),
     $$
     $$
  l_{12}(t_{12})=(\cosh 3Jv/2 - \cosh Jv/2 )/2=\tilde{l}_{12}(\tilde{t}_{12}),  \quad
\tilde{l}_{12}({t}_{21})=ij_2^2J^{-1}\sinh Jv,
     $$
     $$
  \tilde{l}_{12}(t_{12})=-ij_2^2(2J)^{-1}(\sinh 3Jv/2 + \sinh Jv/2 ), \quad
  l_{21}(\tilde{t}_{12})=-ij_1^2J^{-1}\sinh Jv,     \quad
     $$
     $$
  l_{21}(\tilde{t}_{21})=ij_1^2(2J)^{-1}(\sinh 3Jv/2 + \sinh Jv/2 ), \quad
   \tilde{l}_{21}(t_{12})=ij_2^2J^{-1}\sinh Jv,
    $$
    $$
 \tilde{l}_{21}(t_{21})=-ij_2^2(2J)^{-1}(\sinh 3Jv/2 + \sinh Jv/2 ), \quad
 l_{13}(t_{13})=(\cosh 2Jv - 1)/2 =\tilde{l}_{13}(\tilde{t}_{13}),  \quad
    $$
     $$
 l_{21}(t_{21})=(\cosh 3Jv/2 - \cosh Jv/2 )/2=\tilde{l}_{21}(\tilde{t}_{21}),
     $$
     \begin{equation}
 l_{13}(\tilde{t}_{13})=-iJ^{-1}(2\sinh Jv - \sinh 2Jv),             \quad
 \tilde{l}_{13}({t}_{13})=iJ(2\sinh Jv - \sinh 2Jv).
    \label{42}
    \end{equation}
 Only nonzero expressions are written out above.
According to the additional relations (\ref{34}) there are three
independent generators of $ so_{v}(3;j;{\bf C}), $ for example,
     $ l_{11}, l_{12}, \tilde{l}_{12}. $
Their commutation relations follow from Eq.(\ref{33})
$$
l_{11}l_{12} {\cosh}Jv - l_{12} l_{11} =
 l_{11} \tilde l_{12}
ij_1^{2}J^{-1}{\sinh}Jv,
$$
 $$
l_{11}  \tilde  l_{12}  {\cosh}Jv   - \tilde l_{12} l_{11} = -
l_{11} l_{12} ij_2^{2}J^{-1}{\sinh}Jv ,
$$
\begin{equation}
\left [ l_{12} , \tilde l_{12} \right ] =
\left (1 - l_{11}^2 \right )
2iJ{\sinh}Jv/2 -i(j_2^2l_{12}^2 + j_1^2 \tilde l_{12}^2)J^{-1}{\tanh}Jv/2.
\end{equation}

The quantum analogue of the  universal enveloping
 algebra of CK algebra $ so(3;j;{\bf C})=\{X_{01},X_{02},X_{12}\}  $ with
 the rotation generator $ X_{02} $ as the primitive
element of the Hopf algebra has been given in [11],[12].
 The Hopf algebra structure of
 $ so_w(3;j;X_{02}) $  is given by
      $$
\Delta X_{02}  =  I\otimes X_{02}+X_{02}\otimes I,
      $$
      $$
\Delta X  =  e^{-wX_{02}/2}\otimes X+X\otimes e^{wX_{02}/2},
\quad  X=X_{01},X_{12},
      $$
      $$
\epsilon(X_{01})  =  \epsilon(X_{02})=\epsilon(X_{12})=0, \quad
 S(X_{02})  =  -X_{02},
      $$
      $$
 S(X_{01})  =  -X_{01}\cos{Jw/2}+
j_1^2X_{12}J^{-1}\sin{Jw/2},
      $$
      $$
 S(X_{12})  =  -X_{12}\cos{Jw/2}-
j_2^2X_{01}J^{-1}\sin{Jw/2},
      $$
 \begin{equation}
 [X_{01},X_{02}]=j_1^2X_{12},\quad   [X_{02},X_{12}]=j_2^2X_{01},\quad
 [X_{12},X_{01}]={ \sinh{ wX_{02}} \over {w}}.
\label{43}
 \end{equation}

The isomorphism of $ so_w(3;j;X_{02}) $ and quantum algebra
     $ so_v(3;j;{\bf C}) $ is easily established with the help of the
following   relations   between  generators  and  deformation
 parameters
        $$
   l_{11}=e^{-wX_{02}}, \quad
 l_{12}=JEX_{01}e^{-wX_{02}/2}, \quad
 \tilde{l}_{12}= JEX_{12}e^{-wX_{02}/2},
        $$
 \begin{equation}
       v=-iw, \quad
E  =  i{\left ( 2wJ^{-1}\sin Jw \right ) }^{1/2}.
\label{44}
\end{equation}

Now the quantum analogous of the nonsemisimple CK groups and
algebras are obtained by specific values of the parameters
$ j_1,j_2. $ In particular $ j_1=\iota_1, j_2=1 $
corresponds to Euclidean quantum group $ E_{v}(2;{\bf C}) $
(cf. [11]--[13]) and
$ j_1=\iota_1, j_2=\iota_2 $ corresponds to Galilean quantum
group $ G_{v}(2;{\bf C}) $ (cf. [12],[14]).
So the quantum orthogonal CK algebras may be constructed
both as the dual to the quantum group and by the contractions
of quantum orthogonal  algebras.

For $ j_1=j_2=1 $ we have the quantum group
$ SO_q(3;{\bf C}) $ and the quantum algebra
$ so_q(3;{\bf C}). $
Let us mark the elements of the generating
matrix of $ SO_q(3;{\bf C}), $ represented in the
form (\ref{36}), and generators of $ so_q(3;{\bf C}), $
represented in the form (\ref{37}), (\ref{38}), with the
prime. Then all formulas for
$ SO_v(3;j;{\bf C}) $ and $ so_v(3;j;{\bf C}) $
may be obtained from
the corresponding formulas for
$ SO_q(3;{\bf C}) $ and $ so_q(3;{\bf C}) $
by the following transformations of generators and deformation
parameter
$$
t'_{11}=t_{11}, \ {\tilde t}'_{11}=j_1j_2{\tilde t}_{11}, \
t'_{12}=j_1t_{12}, \ {\tilde t}'_{12}=j_2{\tilde t}_{12},
$$
$$
t'_{13}=t_{13}, \ {\tilde t}'_{13}=j_1j_2{\tilde t}_{13}, \
t'_{21}=j_1t_{21}, \ {\tilde t}'_{21}=j_2{\tilde t}_{21},
$$
$$
z=Jv,
$$
$$
l_{11}=l'_{11}, \  l_{12}=j_1 l'_{12}, \
{\tilde l}_{12}=j_2{\tilde l}'_{12}, \ l_{21}=j_1l'_{21}, \
{\tilde l}_{21}=j_2{\tilde l}'_{21},
$$
\begin{equation}
l_{13}= l'_{13}, \
{\tilde l}_{13}=j_1j_2{\tilde l}'_{13}.
\label{53}
\end{equation}
This is nothing else than the contraction transformation if one
replace parameters $ j_k $ with new parameters
$ \epsilon_k, $ which tends to zero.
It worth mention that In\"on\"u-Wigner contractions [9]
of groups at least in CK scheme are just the
regarding of groups over algebras $ {\bf D}_n(j) $
with all or some nilpotent parameters $ j_k. $

\section{Conclusion}

It was demonstrated in previous Sections that the
orthogonal CK groups are naturally arrised as the matrix groups over
the algebra {\bf D} with nilpotent commutative generators.
Their explicit realization depend on the choice of the
basis in the corresponding CK vector space over
the algebra {\bf D}. The Cartesian basis provide the most
simple and well known representations of the CK groups. Then
the realization in an arbitrary basis may be obtained with the
help of the similiraty transformations. In particular,
the realization in  so-called symplectic basis is
needed for quantum deformations of the orthogonal CK groups.
We  have shown that FRT quantization theory of simple
(and semisimple) groups describe also the deformations
of the nonsemisimple orthogonal CK groups and algebras,
if apply it to the corresponding objects over Pimenov algebra {\bf D}).
The ambiguity of the transformations from Cartesian to symplectic
basis provide the different realization of the
diagonal elements of the matrix $ T(j), $ as the
elements of the algebra {\bf D}, which on the level of
quantum algebras leads to the different choice
of the primitive elements of the Hopf algebra.

\section{Acknowledgments}
One of us (NG) would like to thank the Organizing Committee of QGDC
especially Erdal In\"on\"u and Metin Arik for all their help and hospitality.

\end{document}